\documentclass[11pt,a4paper]{article}
\usepackage{amsmath,amsthm,amsfonts,amssymb,mathrsfs}
\usepackage{graphicx}
\usepackage{subfigure}
\usepackage{authblk}
\usepackage[latin1]{inputenc}
\usepackage[swedish,english]{babel}
\newtheorem{remark}{Remark}[section]

\begin{document}
\title{Variants of an explicit kernel-split panel-based Nyström
  discretization scheme for Helmholtz boundary value problems}
\author{Johan Helsing and Anders Holst\\
  {\it Centre for Mathematical Sciences}\\
  {\it Lund University, Box 118, 221 00 Lund, Sweden}}
\date{\today} 
\maketitle

\begin{abstract}
  The incorporation of analytical kernel information is exploited in
  the construction of Nyström discretization schemes for integral
  equations modeling planar Helmholtz boundary value problems.
  Splittings of kernels and matrices, coarse and fine grids,
  high-order polynomial interpolation, product integration performed
  on the fly, and iterative solution are some of the numerical
  techniques used to seek rapid and stable convergence of computed
  fields in the entire computational domain.
\end{abstract}

\section{Introduction}

The question of what high-order accurate Nyström discretization scheme
is the most efficient for solving planar or axisymmetric Laplace,
biharmonic, and Helmholtz boundary value problems, modeled as integral
equations, is a topic of current interest in computational
mathematics. Particularly intriguing are situations where the solution
needs to be evaluated in the entire computational domain, also close
to domain boundaries~\cite{Kloc13}.

The recent paper~\cite{Hao13} classifies Nyström schemes into four
categories depending on whether they are ``global'' or ``panel-based''
and whether they use an ``explicit kernel split'' or ``no explicit
kernel split'' for the discretization of integral operators with
singular kernels. A global Nyström scheme uses the periodic
trapezoidal rule as its underlying quadrature rule in the
discretization. This quadrature has the advantage that exponential
convergence can be obtained provided that certain regularity
assumptions hold on the integrand~\cite[Theorem~12.6]{Kres99}. Global
schemes are therefore the most efficient in many situations.
Panel-based quadrature, such as composite Gauss--Legendre quadrature,
merely achieves polynomial order convergence, but is better suited for
adaptivity and may offer more flexibility in the presence of various
(near) singularities that arise when solution fields are to be
evaluated close to domain boundaries and when domain boundaries are
unions of smooth open arcs. Quadrature schemes that explicitly split
singular kernels into smooth parts and parts with known singularities
may enable higher achievable accuracy and more rapid convergence than
general-purpose schemes which do not use this information.
See~\cite{Hao13} for a general discussion of the merits of different
combinations of discretization strategies and~\cite{Barn14} for recent
progress on explicit kernel-split global Nyström schemes.

The purpose of the present work is to investigate the performance of
explicit kernel-split panel-based Nyström schemes, constructed by
further developing ideas presented
in~\cite{Hels08,Hels09,Hels13a,Hels13b,Hels14}, and to facilitate a
comparison of these new schemes with the split-free panel-based
schemes actually implemented in~\cite{Hao13}. The outcome of a such a
comparison depends, of course, on many things including the test
problem chosen, the details of the implementations, and what aspects
of the schemes that are compared. In the present study we choose to
solve the planar high-frequency exterior Helmholtz Dirichlet problem
of~\cite[Figure 4(c)]{Hao13} and concentrate on convergence speed and
on achievable accuracy in far fields and near fields. While we refrain
from selecting an overall winner, we demonstrate that explicit
kernel-split panel-based schemes are indeed competitive and we provide
a number of numerical tools for enhancing their performance beyond
that of naive implementations.

\section{The exterior Helmholtz Dirichlet problem}

Let $D$ be a bounded simply connected domain in $\mathbb{R}^2$ with
boundary $\gamma$, let $E$ be the exterior to the closure of $D$, let
$r=(x,y)$ be a point of $E\cup\gamma$, and let $\nu$ be the exterior
unit normal to $D$ defined for almost every $r\in\gamma$. The exterior
Dirichlet problem for the Helmholtz equation
\begin{align}
  \Delta u(r)+k^2 u(r)&=0\,,
  \quad r\in E\,,\label{eq:Helm}\\
\lim_{E\ni r\to r^{\circ}}
  u(r)&=g(r^{\circ})\,,\quad r^{\circ}\in\gamma\,,\label{eq:HelmBV}\\
  \lim_{|r|\to\infty}
  \sqrt{|r|}\left(\frac{\partial}{\partial|r|}-{\rm
        i}k\right)u(r)&=0\,,
\end{align}
has a unique solution $u(r)$ under mild assumptions on $\gamma$ and
$g(r)$~\cite{Mitr96} and can be modeled using a combined field
integral representation~\cite[Equation~(3.25)]{Colt98} in terms of a
layer density $\rho(r)$
\begin{equation}
u(r)=\int_{\gamma}\frac{\partial\Phi_k}{\partial\nu'}(r,r')\rho(r')
\,{\rm d}\sigma'-
\frac{{\rm i}k}{2}\int_{\gamma}\Phi_k(r,r')\rho(r')
\,{\rm d}\sigma'\,, \quad r\in E\,.
\label{eq:Helmrep1}
\end{equation}
Here ${\rm d}\sigma$ is an element of arc length, differentiation with
respect to $\nu$ denotes the normal derivative, and $\Phi_k(r,r')$ is
the fundamental solution to~(\ref{eq:Helm}) 
\begin{equation}
\Phi_k(r,r')=\frac{\rm i}{4}H_0^{(1)}(k|r-r'|)\,,
\end{equation}
where $H_0^{(1)}$ is the zeroth order Hankel function of the first kind.

Insertion of~(\ref{eq:Helmrep1}) into~(\ref{eq:HelmBV}) gives the
combined field integral equation
\begin{equation}
\left(I+K_k-\frac{{\rm i}k}{2}S_k\right)\rho(r)=2g(r)\,, 
\quad r\in\gamma\,,
\label{eq:HelmDiri1}
\end{equation}
where
\begin{align}
  K_k\rho(r)&=
  2\int_{\gamma}\frac{\partial\Phi_k}{\partial\nu'}(r,r')\rho(r')
  \,{\rm d}\sigma'\,,\label{eq:Koper}\\
  S_k\rho(r)&=2\int_{\gamma}\Phi_k(r,r')\rho(r')\,
  {\rm d}\sigma'\,.
\end{align}

\begin{remark}
  {\rm The representation~(\ref{eq:Helmrep1}) contains a real valued
    coupling parameter, denoted $\eta$
    in~\cite[Equation~(3.25)]{Colt98}, which we have set to
    $\eta=k/2$. The choice of $\eta$ may greatly influence the
    spectral properties of $I+K_k-{\rm i}\eta S_k$ and affect the
    achievable accuracy in solutions to discretized versions
    of~(\ref{eq:HelmDiri1}). Convergence rates of iterative solvers
    are affected, too. See~\cite[Section IIB]{Betc11} for a review of
    recommendations for $\eta$ when $D$ is a starlike domain.}
\label{rmk:eta}
\end{remark}

\section{Panel-based Nyström discretization}
\label{sec:pan}

Let us think of $K_k-{\rm i}kS_k/2$ as a single integral operator $M$
with kernel $M(r,r')$ and add subscript ``$\gamma$'' or ``$E$'' when
it is instructive to point out if $r\in\gamma$ or $r\in E$.
Equations~(\ref{eq:HelmDiri1}) and~(\ref{eq:Helmrep1}) then assume the
general form
\begin{gather}
  \rho(r)+\int_{\gamma}M_{\gamma}(r,r')\rho(r')\,{\rm d}\sigma'=2g(r)\,,
  \quad r\in\gamma\,,\label{eq:HD}\\
  u(r)=\frac{1}{2}\int_{\gamma}M_E(r,r')\rho(r')\,{\rm d}\sigma'\,,
  \quad r\in E\,.\label{eq:Hr}
\end{gather}

An $n_{\rm pt}$-point panel-based Nyström discretization scheme
for~(\ref{eq:HD}) and~(\ref{eq:Hr}) involves the following setup
steps: choose a parameterization $r(t)$ of $\gamma$; construct a mesh
of $n_{\rm pan}$ quadrature panels on $\gamma$; choose an underlying
interpolatory quadrature rule with nodes $\mathfrak{t}_i$ and weights
$\mathfrak{w}_i$, $i=1,\ldots,n_{\rm pt}$, on a canonical interval
$[-1,1]$; find actual nodes $t_i$ and weights $w_i$,
$i=1,\ldots,n_{\rm pt}n_{\rm pan}$, on the panels of $\gamma$ via
transformations of $\mathfrak{t}_i$ and $\mathfrak{w}_i$.

The discretization scheme could proceed with an approximation of the
integrals in~(\ref{eq:HD}) using $t_i$ and $w_i$, and the demand that
the discretization holds at the nodes $t_i$. Introducing the speed
function $s(t)=|\dot{r}|=|{\rm d}r(t)/{\rm d}t|$, the resulting system
would look like
\begin{equation}
\rho_i+\sum_{j=1}^nM_{\gamma}(r_i,r_j)\rho_js_jw_j=2g_i\,,\quad i=1,\ldots,n\,,
\label{eq:HD1}
\end{equation}
where $n=n_{\rm pt}n_{\rm pan}$, $r_i=r(t_i)$, $\rho_i=\rho(r(t_i))$,
$s_i=s(t_i)$, and $g_i=g(r(t_i))$. Upon solving~(\ref{eq:HD1}) for
$\rho_i$, the field $u(r)$ could be obtained from a discretization
of~(\ref{eq:Hr})
\begin{equation}
u(r)=\frac{1}{2}\sum_{j=1}^nM_E(r,r_j)\rho_js_jw_j\,.
\label{eq:Hr1}
\end{equation}
Note that the grid points $r(t_i)$ in~(\ref{eq:HD1}) both play the
role of {\it target points} $r_i$ and of {\it source points} $r_j$.

The simple scheme of~(\ref{eq:HD1}) and~(\ref{eq:Hr1}) works well if
$\gamma$, $M(r,r')$, and $g(r)$ are smooth. Then
$M_{\gamma}(r_i,r(t))$ and $M_E(r,r(t))$ are well approximated by
polynomials in $t$. One can show that under suitable regularity
assumptions on $M(r,r')$ and $\rho(r)$, the convergence rate of
Nyström schemes reflect those of their underlying
quadratures~\cite[Section~12.2]{Kres99}. An underlying $n_{\rm
  pt}$-point Gauss--Legendre quadrature would result in a scheme of
order $2n_{\rm pt}$.

Now, for Helmholtz problems, $M(r,r')$ is not smooth. Depending on how
$r'$ approaches $r$, the kernels of $K_k$ and $S_k$ can contain both
logarithmic- and Cauchy-type singularities. Since such singularities
are difficult to resolve by polynomials, the convergence of the
scheme~(\ref{eq:HD1}) and~(\ref{eq:Hr1}) will be slow. In the context
of panel-based schemes it is therefore common to single out quadrature
panels where some special-purpose quadrature is required for
efficiency. The scheme~(\ref{eq:HD1}) and~(\ref{eq:Hr1}) then assumes
the form
\begin{gather}
\begin{split}
\rho_i+\sum_{j\in{\cal C}(r_i)}M_{\gamma}&(r_i,r_j)\rho_js_jw_{ij}\\
      +&\sum_{j\in{\cal F}(r_i)}M_{\gamma}(r_i,r_j)\rho_js_jw_j
=2g_i\,,\quad i=1,\ldots,n\,,
\end{split}
\label{eq:HD2}\\
u(r)=\frac{1}{2}\sum_{j\in{\cal C}(r)}M_E(r,r_j)\rho_js_jw_{ij}
    +\frac{1}{2}\sum_{j\in{\cal F}(r)}M_E(r,r_j)\rho_js_jw_j\,.
\label{eq:Hr2}
\end{gather}
Here ${\cal C}(r_i)$ and ${\cal C}(r)$ are sets of source points on
panels that are close to $r_i$ and $r$, respectively, and where
special-purpose quadrature weights $w_{ij}$ are used. Source points on
remaining panels are contained in the sets ${\cal F}(r_i)$ and ${\cal
  F}(r)$. These panels are considered to be sufficiently far away from
$r_i$ and $r$ for the kernels to be smooth and for the underlying
quadrature to be efficient. In the present work, the set ${\cal
  C}(r_i)$ contains source points on at most three panels: the panel
on which $r_i$ is situated and one or both of its neighboring panels.

See~\cite{Hao13} for a review of special-purpose quadratures that can
be used for $M_{\gamma}$. We intend to use product integration derived
with polynomial interpolation. In general, strictly panel-based
product integration of this type does not result in more than order
$n_{\rm pt}$ convergence~\cite[Section 4.2.2]{Atki97}.

\section{The partitioning and splitting of matrices}

This section casts the linear system~(\ref{eq:HD2}) and the
post-processor~(\ref{eq:Hr2}) into matrix-vector form and introduces
matrix splittings that help to simplify the description of our
discretization schemes.

The system~(\ref{eq:HD2}) can be written
\begin{equation}
\left({\bf I}+{\bf M}_{\gamma}\right)\boldsymbol{\rho}=2{\bf g}\,,
\label{eq:HD3}
\end{equation}
where ${\bf I}$ is the $n\times n$ identity matrix, ${\bf M}_{\gamma}$
is an $n\times n$ matrix containing the discretization of
$M_{\gamma}$, and $\boldsymbol{\rho}$ and ${\bf g}$ are column vectors
of length $n$ containing discrete values of $\rho(r)$ and $g(r)$.
Based on panel affiliation one can partition ${\bf M}_{\gamma}$ into
$n_{\rm pan}\times n_{\rm pan}$ square blocks with $n_{\rm pt}^2$
entries each. One can also split $M_{\gamma}(r,r')$ into two functions
\begin{equation}
M_{\gamma}(r,r')=M_{\gamma}^{\star}(r,r')+M_{\gamma}^{\circ}(r,r')\,.
\label{eq:Msplit}
\end{equation}
Here $M_{\gamma}^{\star}(r,r')$ is zero except for when $r$ and $r'$
are situated on the same or on neighboring panels. In this latter case
$M_{\gamma}^{\circ}(r,r')$ is zero. The kernel
splitting~(\ref{eq:Msplit}) corresponds to a matrix splitting and we
can write~(\ref{eq:HD3}) as
\begin{equation}
\left({\bf I}+{\bf M}_{\gamma}^{\star}
      +{\bf M}_{\gamma}^{\circ}\right)\boldsymbol{\rho}=2{\bf g}\,,
\label{eq:HD4}
\end{equation}
where ${\bf M}_{\gamma}^{\star}$ is the block-tridiagonal part of
${\bf M}_{\gamma}$ plus its upper right and lower left blocks. Note
that ${\bf M}_{\gamma}^{\star}$ contains all entries of ${\bf
  M}_{\gamma}$ which involve special-purpose quadrature (product
integration) along with some other entries for which the underlying
quadrature is sufficient. Since ${\bf M}_{\gamma}^{\star}$ only has
$3n_{\rm pt}n$ non-zero entries one can say that~(\ref{eq:HD4}) is of
FMM-compatible form~\cite[Definition~1.1]{Hao13}.

Assuming that the field $u(r)$ is to be evaluated at $n_{\rm fp}$
different field points $r\in E$ one can write~(\ref{eq:Hr2}) in the
form
\begin{equation}
{\bf u}=\frac{1}{2}{\bf M}_E\boldsymbol{\rho}\,,
\label{eq:Hr3}
\end{equation}
where ${\bf u}$ is a column vector with $n_{\rm fp}$ entries and ${\bf
  M}_E$ is a $n_{\rm fp}\times n$ rectangular matrix. Collecting the
entries of ${\bf M}_E$ that correspond to the first sum
in~(\ref{eq:Hr2}) in a matrix ${\bf M}_E^{\star}$ and the remaining
entries in a matrix ${\bf M}_E^{\circ}$ we write~(\ref{eq:Hr3}) as
\begin{equation}
{\bf u}=\frac{1}{2}\left({\bf M}_E^{\star}+{\bf M}_E^{\circ}\right)
\boldsymbol{\rho}\,.
\label{eq:Hr4}
\end{equation}

If $n_{\rm pt}$-point Gauss--Legendre quadrature is used as the
underlying quadrature, then the actions of ${\bf M}_{\gamma}^{\circ}$
and ${\bf M}_E^{\circ}$ on $\boldsymbol{\rho}$ in~(\ref{eq:HD4})
and~(\ref{eq:Hr4}) correspond to $2n_{\rm pt}$th order accurate
discretization of distant interactions in $M$. The action of ${\bf
  M}_{\gamma}^{\star}$ and ${\bf M}_E^{\star}$ on $\boldsymbol{\rho}$
corresponds to an, at most, $n_{\rm pt}$th order accurate
discretization of close interactions in $M$ and limits the overall
convergence rate of the Nyström scheme. It is therefore important to
make the asymptotic error constant of this discretization small.

\section{The known singularities in $S_k$ and $K_k$}
\label{sec:known}

This section reviews singularities that arise in the kernels of the
operators $S_k$ and $K_k$ as $r'\in\gamma$ approaches $r$. Knowledge
of these singularities is essential for constructing explicit-split
discretization schemes. A similar review can be found
in~\cite[Section~3.5]{Colt98}.

The kernel of $S_k$ can be expressed in the form
\begin{equation}
S_k(r,r')=S_{0k}(r,r')-\frac{2}{\pi}
\log|r-r'|\Im\left\{S_k(r,r')\right\}\,,
\label{eq:Sform}
\end{equation}
where $S_{0k}(r,r')$ and $\Im\left\{S_k(r,r')\right\}$ are smooth
functions with limits
\begin{gather}
\lim_{r'\to r}S_{0k}(r,r')=\frac{\rm i}{2}
-\frac{1}{\pi}\left(\log\left|\frac{k}{2}\right|-\psi(1)\right)\,,\\
\lim_{r'\to r}\Im\left\{S_k(r,r')\right\}=\frac{1}{2}\,.
\end{gather}
Here $\psi$ is the digamma function.

The kernel of $K_k$ can, for $r\in\gamma$, be expressed in a form
analogous to~(\ref{eq:Sform})
\begin{equation}
K_k(r,r')=K_{0k}(r,r')-\frac{2}{\pi}
\log|r-r'|\Im\left\{K_k(r,r')\right\}\,,
\label{eq:Kform1}
\end{equation}
where $K_{0k}(r,r')$ and $\Im\left\{K_k(r,r')\right\}$ are smooth
functions with limits
\begin{gather}
  \lim_{r'\to r}K_{0k}(r,r')=\frac{1}{2\pi}
  \frac{(\nu\cdot\ddot{r})}{|\dot{r}|^2}\,,\\
  \lim_{r'\to r}\Im\left\{K_k(r,r')\right\}=0\,.
\end{gather}
In~(\ref{eq:Kform1}) we use $\dot{r}={\rm d}r(t)/{\rm d}t$ and
$\ddot{r}={\rm d}^2r(t)/{\rm d}t^2$.

The kernel of $K_k$ can, for $r\in E$, be expressed in the form
\begin{equation}
K_k(r,r')=K_{0k}(r,r')-\frac{2}{\pi}
\log|r-r'|\Im\left\{K_k(r,r')\right\}
-\frac{1}{\pi}\frac{((r'-r)\cdot\nu')}{|r'-r|^2}\,,
\label{eq:Kform2}
\end{equation}
where $K_{0k}(r,r')$ and $\Im\left\{K_k(r,r')\right\}$ are smooth
functions. The expression~(\ref{eq:Kform2}), which was derived
in~\cite[Section III{\it B}]{Hels13b}, can be verified via the
definition of $K_k$ in~(\ref{eq:Koper}) and a series representation of
the first order Hankel function of the first kind
$H^{(1)}_1(k|r-r'|)$.

\section{Product integration for singular integrals}
\label{sec:prod}

This section reviews a special-purpose quadrature applicable to the
singular kernels of $S_k$ and $K_k$. The presentation is a summary
of~\cite[Section~9]{Hels14} and concerns the discretization of the
integral
\begin{equation}
I_p(r)=\int_{\gamma_p}G(r,r')\rho(r')\,{\rm d}\sigma'\,,
\label{eq:ex1}
\end{equation}
where $G(r,r')$ is a non-smooth kernel, $\rho(r)$ is a smooth layer
density, $\gamma_p$ is a quadrature panel on a curve $\gamma$ with
endpoints $r(t_a)$ and $r(t_b)$, $t_a<t_b$, and the target point $r$
is located close to, or on, $\gamma_p$. Gauss--Legendre quadrature is
used as underlying quadrature with nodes $t_i\in[t_a,t_b]$ and weights
$w_i$, $i=1,\ldots,n_{\rm pt}$. For brevity we write
$\rho(t)=\rho(r(t))$.

\subsection{Logarithmic singularity plus smooth part}

Consider~(\ref{eq:ex1}) when $G(r,r')$ can be expressed as
\begin{equation}
G(r,r')=G_0(r,r')+\log|r-r'|G_{\rm L}(r,r')\,,
\label{eq:ex2}
\end{equation}
where both $G_0(r,r')$ and $G_{\rm L}(r,r')$ are smooth functions.
Then one can find, using polynomial product integration against the
logarithmic kernel~\cite[Section~2.3]{Hels09}, weight corrections
$w_{{\rm L}j}^{\rm corr}(r)$ such that
\begin{equation}
I_p(r)=\sum_{j=1}^{n_{\rm pt}} G(r,r_j)\rho_js_jw_j
+\sum_{j=1}^{n_{\rm pt}} G_{\rm L}(r,r_j)\rho_js_jw_jw_{{\rm L}j}^{\rm corr}(r)
\label{eq:fix2}
\end{equation}
is exact for $G_0(r,r(t))\rho(t)s(t)$ being a polynomial of degree
$2n_{\rm pt}-1$ in $t$ and for $G_{\rm L}(r,r(t))\rho(t)s(t)$ being a
polynomial of degree $n_{\rm pt}-1$.

While it is rather easy to compute $w_{{\rm L}j}^{\rm corr}(r)$ for a
general point $r$, it its even easier in the special case that $r$
coincides with a target point $r_i$ on $\gamma_p$.
Then~(\ref{eq:fix2}) becomes
\begin{equation}
I_p(r_i)=\sum_{j\ne i}^{n_{\rm pt}}G(r_i,r_j)\rho_js_jw_j
+G_0(r_i,r_i)\rho_is_iw_i
+\sum_{j=1}^{n_{\rm pt}}
G_{\rm L}(r_i,r_j)\rho_js_jw_jw_{{\rm L}j}^{\rm corr}(r_i)\,,
\end{equation}
where
\begin{equation}
w_{{\rm L}j}^{\rm corr}(r_i)=\left\{
\begin{array}{lr}
\mathfrak{W}_{{\rm L}ij}/\mathfrak{w}_j
-\log\left|\mathfrak{t}_i-\mathfrak{t}_j\right|\,, & j\ne i\,,\\
\mathfrak{W}_{{\rm L}ii}/\mathfrak{w}_i
+\log\left|(t_b-t_a)s_i/2\right|\,, & j=i\,.
\end{array}
\right.
\label{eq:Acorr}
\end{equation}
Here $\mathfrak{W}_{\rm L}$ is a square matrix whose entries are
$(n_{\rm pt}-1)$th degree product integration weights for the
logarithmic integral operator on the canonical interval and only
depend on the nodes $\mathfrak{t}_i$. See Section~\ref{sec:pan} for
definitions of $\mathfrak{t}_i$ and $\mathfrak{w}_i$.

Note that the off-diagonal corrections in~(\ref{eq:Acorr}) do not
depend on $\gamma_p$ and that $\mathfrak{W}_{\rm L}$ only needs to be
computed and stored once. An analogous derivation for $r_i$ and $r_j$
on neighboring panels shows that the corresponding corrections then
depend on the nodes $\mathfrak{t}_i$ and the relative length of the
panels. Appendix~A contains a {\sc Matlab} function that constructs
the matrix $\mathfrak{W}_{\rm L}$.

\subsection{Logarithmic- and Cauchy-type singularities plus smooth part}

Now consider~(\ref{eq:ex1}) when $G(r,r')$ can be expressed as
\begin{equation}
G(r,r')=G_0(r,r')+\log|r-r'|G_{\rm L}(r,r')
+\frac{(r'-r)\cdot\nu'}{|r'-r|^2}G_{\rm C}(r,r')\,,
\label{eq:ex3}
\end{equation}
where $G_0(r,r')$, $G_{\rm L}(r,r')$, and $G_{\rm C}(r,r')$ are smooth
functions. If $r\in\gamma$, then the third term on the right
in~(\ref{eq:ex3}) is a smooth function and we are back
to~(\ref{eq:ex2}). Otherwise one can find, using polynomial product
integration against the Cauchy-singular
kernel~\cite[Section~2.1]{Hels09}, compensation weights $w_{{\rm
    C}j}^{\rm cmp}(r)$ such that
\begin{equation}
I_p(r)=\sum_{j=1}^{n_{\rm pt}} G(r,r_j)\rho_js_jw_j
+\sum_{j=1}^{n_{\rm pt}}G_{\rm L}(r,r_j)\rho_js_jw_jw_{{\rm L}j}^{\rm corr}(r)
+\sum_{j=1}^{n_{\rm pt}}G_{\rm C}(r,r_j)\rho_jw_{{\rm C}j}^{\rm cmp}(r)
\label{eq:fix3}
\end{equation}
is exact under the same conditions as~(\ref{eq:fix2}) and the
additional condition that $G_{\rm C}(r,r(t))\rho(t)s(t)$ is a
polynomial of degree $n_{\rm pt}-1$. Appendix~B contains a {\sc
  Matlab} function that constructs $w_{{\rm L}j}^{\rm corr}(r)$ and
$w_{{\rm C}j}^{\rm cmp}(r)$ for $r\in E$:

\subsection{When to activate product integration}
\label{sec:active}

Due to its comparably low order, special-purpose quadrature for source
points on a panel $\gamma_p$ with endpoints $r(t_a)$ and $r(t_b)$ and
arc length $|\gamma_p|$ should only be activated when it is expected
to give better accuracy than the underlying quadrature. In the
numerical examples of Section~\ref{sec:numex} we use either $n_{\rm
  pt}=16$ or $n_{\rm pt}=32$. For target points $r(t_i)\in\gamma$,
product integration is activated when
\begin{itemize}
\item[$\bullet$] $n_{\rm pt}=16$ and $|t_i-(t_a+t_b)/2|<t_b-t_a$,
\item[$\bullet$] $n_{\rm pt}=32$ and $|t_i-(t_a+t_b)/2|<0.7(t_b-t_a)$.
\end{itemize}
For field points $r\in E$, product integration is activated when
\begin{itemize}
\item[$\bullet$] $n_{\rm pt}=16$ and the minimum distance from $r$ to
  $\gamma_p$ is less than $1.1|\gamma_p|$.
\item[$\bullet$] $n_{\rm pt}=32$ and the minimum distance from $r$ to
  $\gamma_p$ is less than $0.3|\gamma_p|$.
\end{itemize}

\section{Resolution, grids, and interpolation matrices}

The accurate discretization of~(\ref{eq:HD}) and~(\ref{eq:Hr})
requires that the integrand $M(r,r')\rho(r')$ is resolved with respect
to the variable of integration. This task, on a single grid, is often
more expensive than the task of resolving $M(r,r')$ and $\rho(r')$
separately on different grids. There is an obvious resolution level
needed for an accurate discrete representation of $\rho(r')$ while the
resolution level needed for $M(r,r')$ varies with $r-r'$. It would
therefore be beneficial if $M(r,r')$ and $\rho(r')$, somehow, could be
decoupled early in the discretization process -- prior to invoking a
linear solver. Simple tools for achieving this are now provided.
See~\cite{Ho12} for far more advanced schemes based on multilevel
matrix compression.

\subsection{Two grids on $\gamma$}

The product integration of Section~\ref{sec:prod} incorporates
analytical information about the nature of the singularities in
$G(r,r')$ but disregards analytical information about the multiplying
functions $G_{\rm L}(r,r')$ and $G_{\rm C}(r,r')$. A central theme in
the present work is the incorporation of available information into
discretization schemes and we shall not overlook the information
contained in multiplying functions when dealing with $M(r,r')$. Rather
than using this information analytically, however, we exploit it
numerically via a two-grid procedure where, loosely speaking, a coarse
grid is used to resolve $\rho(r')$ and $M(r,r')$ in most situations
but a fine grid is used when $r'$ is close to $r$.
Compare~\cite[Section~5.1]{Hao13}, where several fine grids of
auxiliary points are used for the same purpose.

The construction of our two grids and their accompanying quadratures
is simple, given a mesh of $n_{\rm pan}$ panels on $\gamma$. Nodes
$t_i^{(1)}$, weights $w_i^{(1)}$, and points $r_i^{(1)}$ of the coarse
grid are identical to the quantities $t_i$, $w_i$, and $r_i$
constructed in Section~\ref{sec:pan}. Nodes $t_i^{(2)}$, weights
$w_i^{(2)}$ and points $r_i^{(2)}$ of the fine grid are obtained in an
analogous fashion, but with twice the number of nodes
$\mathfrak{t}_i^{(2)}$ and weights $\mathfrak{w}_i^{(2)}$,
$i=1,\ldots,2n_{\rm pt}$, on the canonical interval.

\subsection{Matrices for panelwise interpolation}
\label{sec:three}

We need discrete operators that perform polynomial interpolation
between functions on the two grids. For this, we introduce the
Vandermonde matrices ${\bf V}^{(11)}$, ${\bf V}^{(12)}$, ${\bf
  V}^{(21)}$, and ${\bf V}^{(22)}$ with entries
\begin{align}
V^{(11)}_{ij}&=\left(\mathfrak{t}_i^{(1)}\right)^{j-1}\,,\quad
i,j=1,\ldots,n_{\rm pt}\,,\\
V^{(12)}_{ij}&=\left(\mathfrak{t}_i^{(1)}\right)^{j-1}\,,\quad
i=1,\ldots,n_{\rm pt}\,,\quad j=1,\ldots,2n_{\rm pt}\,, \\
V^{(21)}_{ij}&=\left(\mathfrak{t}_i^{(2)}\right)^{j-1}\,,\quad
i=1,\ldots,2n_{\rm pt}\,,\quad j=1,\ldots,n_{\rm pt}\,, \\
V^{(22)}_{ij}&=\left(\mathfrak{t}_i^{(2)}\right)^{j-1}\,,\quad
i,j=1,\ldots,2n_{\rm pt}\,.
\end{align}
We then construct the rectangular matrices
\begin{align}
{\bf P}^{(21)}&={\bf V}^{(21)}\left({\bf V}^{(11)}\right)^{-1}\,,
\label{eq:P21}\\
{\bf Q}^{(12)}&={\bf V}^{(12)}\left({\bf V}^{(22)}\right)^{-1}\,,
\label{eq:Q12}
\end{align}
and expand them into rectangular block diagonal matrices ${\bf P}$ and
${\bf Q}$ by $n_{\rm pan}$ times replicating ${\bf P}^{(21)}$ and
${\bf Q}^{(12)}$. Using {\sc Matlab}-style notation this expansion can
be expressed as
\begin{align}
{\bf P}&=
{\tt blkdiag}({\bf P}^{(21)},{\bf P}^{(21)},\ldots,{\bf P}^{(21)})\,,
\label{eq:P}\\
{\bf Q}&=
{\tt blkdiag}({\bf Q}^{(12)},{\bf Q}^{(12)},\ldots,{\bf Q}^{(12)})\,,
\label{eq:Q}
\end{align}

The matrices ${\bf P}$ and ${\bf Q}$ are simple to interpret. When
${\bf P}$ acts from the left on a column vector it performs panelwise
$(n_{\rm pt}-1)$-degree polynomial interpolation from the coarse grid
to the fine grid. In the context of a Nyström method based on Gaussian
quadrature, this could lead to loss of information. Assume, for
example, that a column vector $\boldsymbol{\rho}^{(1)}$ contains
$2n_{\rm pt}$th order accurate entries. Then the entries of ${\bf
  P}\boldsymbol{\rho}^{(1)}$ are only $n_{\rm pt}$th order accurate.
When ${\bf Q}$ acts from the left on a column vector it performs
panelwise $(2n_{\rm pt}-1)$-degree interpolation from the fine grid to
the coarse grid. If $\boldsymbol{\rho}^{(2)}$ contains $2n_{\rm pt}$th
order accurate entries, then the accuracy in ${\bf
  Q}\boldsymbol{\rho}^{(2)}$ is retained. 

\begin{remark}
  The condition numbers of the Vandermonde matrices, needed for the
  construction of ${\bf P}$ and ${\bf Q}$, are high. Still, very
  accurate interpolation can be obtained by ${\bf P}$ and ${\bf Q}$ if
  explicit inverses are avoided and a backward stable solver ({\sc
    Matlab}'s backslash) is used in~(\ref{eq:P21}) and~(\ref{eq:Q12}).
  See~\cite[Appendix~A]{Hels08}.
\end{remark}

\subsection{Extended interpolation}
\label{sec:ext}

Let $\gamma_{p-1}$, $\gamma_p$, and $\gamma_{p+1}$ be three
consecutive quadrature panels on $\gamma$ with endpoints $r(t_a)$,
$r(t_b)$, $r(t_c)$, $r(t_d)$, and $t_a<t_b<t_c<t_d$. Let $n_{\rm s}$
be a small integer and define the extended set of $n_{\rm pt}+2n_{\rm
  s}$ nodes
\begin{equation}
\mathfrak{t}_i^{(1{\rm x})}=\left\{
\begin{array}{ll}
\alpha(\mathfrak{t}_{i+n_{\rm pt}-n_{\rm s}}^{(1)}-1)-1\,, 
& i=1,\ldots,n_{\rm s}\,,\\
\mathfrak{t}_{i-n_{\rm s}}^{(1)}\,, 
& i=n_{\rm s}+1,\ldots,n_{\rm s}+n_{\rm pt}\,,\\
\beta(\mathfrak{t}_{i-n_{\rm pt}-n_{\rm s}}^{(1)}+1)+1\,, 
& i=n_{\rm s}+n_{\rm pt}+1,\ldots,n_{\rm pt}+2n_{\rm s}\,,
\end{array}
\right.
\end{equation}
where $\alpha=(t_b-t_a)/(t_c-t_b)$ and $\beta=(t_d-t_c)/(t_c-t_b)$.

One can now construct the extended Vandermonde matrices ${\bf
  V}_p^{(11{\rm x})}$ and ${\bf V}_p^{(21{\rm x})}$
\begin{align}
V^{(11{\rm x})}_{pij}&=\left(\mathfrak{t}_i^{(1{\rm x})}\right)^{j-1}\,,\quad
i,j=1,\ldots,n_{\rm pt}+2n_{\rm s}\,,\\
V^{(21{\rm x})}_{pij}&=\left(\mathfrak{t}_i^{(2)}\right)^{j-1}\,,\quad
i=1,\ldots,2n_{\rm pt}\,,\quad j=1,\ldots,n_{\rm pt}+2n_{\rm s}\,,
\end{align}
and the extended matrix
\begin{equation}
{\bf P}_p^{(21{\rm x})}={\bf V}_p^{(21{\rm x})}\left({\bf V}_p^{(11{\rm x})}\right)^{-1}\,.
\end{equation}
Note that the matrix ${\bf P}_p^{(21{\rm x})}$ depends on $p$, via
$\alpha$ and $\beta$, whenever adjacent quadrature panels differ in
parameter length.

If one replaces the diagonal blocks ${\bf P}^{(21)}$ and some
neighboring zeros in ${\bf P}$ of~(\ref{eq:P}) with the corresponding,
slightly wider, blocks ${\bf P}_p^{(21{\rm x})}$ one gets a matrix
${\bf P}_{\rm x}$ which, when acting from the left on a column vector
$\boldsymbol{\rho}^{(1)}$, performs $(n_{\rm pt}+2n_{\rm s}-1)$-degree
polynomial interpolation to the fine grid. The interpolated values
${\bf P}_{\rm x}\boldsymbol{\rho}^{(1)}$ on a given panel $\gamma_p$
are determined by $n_{\rm pt}$ values of $\boldsymbol{\rho}^{(1)}$ on
$\gamma_p$ and by an additional $2n_{\rm s}$ values of
$\boldsymbol{\rho}^{(1)}$ on the neighboring panels $\gamma_{p-1}$ and
$\gamma_{p+1}$, so the interpolation is not strictly panelwise.

\section{Four schemes}
\label{sec:four}

Equipped with underlying and special-purpose quadrature, matrix
splittings, criteria for quadrature activation, coarse and fine grids,
and interpolation matrices, we are now in a position to present
meaningful discretization schemes for~(\ref{eq:HD}) and~(\ref{eq:Hr}).
In doing so we indicate coarse and fine grids with superscripts
``(1)'' and ``(2)'', respectively, and points $r\in E$ with
superscript ``(3)''. Discretized integral operators have two
superscripts where the first refers to their points of evaluation and
the second to their source points.

\subsection{Scheme A}

A simple recipe for the discretization of~(\ref{eq:HD})
and~(\ref{eq:Hr}) is to let the mesh on $\gamma$ have panels of equal
length in parameter, take~(\ref{eq:HD3}) and~(\ref{eq:Hr3}) as they
stand, and only use the discretization points of the coarse grid 
\begin{gather}
\left({\bf I}^{(11)}+{\bf M}^{(11)}_{\gamma}\right)\boldsymbol{\rho}^{(1)}
=2{\bf g}^{(1)}\,,\label{eq:A1}\\
{\bf u}^{(3)}=\frac{1}{2}{\bf M}_E^{(31)}\boldsymbol{\rho}^{(1)}
\label{eq:A2}\,.
\end{gather}

\subsection{Scheme B}

Better resolution of $M(r,r')$ when $r'$ is close to $r$ does not
improve the convergence order, but should decrease the error constant
\begin{gather}
\left({\bf I}^{(11)}
+{\bf Q}{\bf M}^{\star(22)}_{\gamma}{\bf P}
+{\bf M}^{\circ(11)}_{\gamma}
\right)\boldsymbol{\rho}^{(1)}
=2{\bf g}^{(1)}\,,\label{eq:B1}\\
{\bf u}^{(3)}=\frac{1}{2}\left(
 {\bf M}_E^{\star(32)}{\bf P}
+{\bf M}_E^{\star\circ(32)}{\bf P}
+{\bf M}_E^{\circ(31)}\right)\boldsymbol{\rho}^{(1)}
\label{eq:B2}\,.
\end{gather}
Here ${\bf M}_E^{\star\circ(32)}$ is a matrix whose non-zero entries
correspond to interaction between points $r\in E$ and source points on
panels $\gamma_p$ not accounted for in ${\bf M}_E^{\star(32)}$ or
${\bf M}_E^{\circ(31)}$. Scheme B is the one of our schemes that most
resembles the split-free panel-based scheme called ``Modified
Gaussian'' in~\cite{Hao13} and which uses Kolm--Rokhlin special-purpose
quadrature~\cite{Kolm01}.

\subsection{Scheme C}

Scheme C is the same as Scheme B, but with panels that are equal in
arc length $\sigma$ and a with unit speed parameterization $r_{\rm
  uni}(\sigma)$ of $\gamma$. Given any parameterization $r(t)$ of
$\gamma$, it is easy to construct a mesh with panels that are equal in
parameter length $t$. An advantage with such a mesh is that the number
of different special-purpose weights $w_{ij}$, needed
in~(\ref{eq:HD2}) and~(\ref{eq:Hr2}), is low. Still, unless $r(t)$ has
unit speed, equal parameter length panels are not equal in arc length
and this leads to an (unwanted) difference in spacing between
discretization points $r_i$ on different parts of $\gamma$ which, in
turn, may delay convergence. Fortunately, it is not necessary to have
access to $r_{\rm uni}(\sigma)$ in closed form in order to find $r_i$
consistent with a unit speed parameterization. Given any reasonable
parameterization $r(t)$ in closed form, values of the function
$\sigma(t)$ and its inverse can be computed numerically to machine
precision using quadrature and Newton's method. Thus, it is simple to
find parameter values $t_i$ corresponding to desired nodes $\sigma_i$
and $r_i=r_{\rm uni}(\sigma_i)=r(t_i)$. Derivatives of $r_{\rm
  uni}(\sigma)$ with respect to $\sigma$ can be computed from $r(t)$
using elementary calculus.

\subsection{Scheme D}

Replacement of ${\bf P}$ in~(\ref{eq:B1}) and~(\ref{eq:B2}) with the
extended interpolation matrix ${\bf P}_{\rm x}$ of
Section~\ref{sec:ext} results in the $(n_{\rm pt}+2n_{\rm s})$th order
accurate scheme
\begin{gather}
\left({\bf I}^{(11)}
+{\bf Q}{\bf M}^{\star(22)}_{\gamma}{\bf P}_{\rm x}
+{\bf M}^{\circ(11)}_{\gamma}
\right)\boldsymbol{\rho}^{(1)}
=2{\bf g}^{(1)}\,,\label{eq:D1}\\
{\bf u}^{(3)}=\frac{1}{2}\left(
 {\bf M}_E^{\star(32)}{\bf P}_{\rm x}
+{\bf M}_E^{\star\circ(32)}{\bf P}_{\rm x}
+{\bf M}_E^{\circ(31)}\right)\boldsymbol{\rho}^{(1)}
\label{eq:D2}\,.
\end{gather}
Scheme D is implemented with equal arc length panels, a unit speed
parameterization of $\gamma$, and $n_{\rm s}=4$.

\section{Numerical examples}
\label{sec:numex}

This section tests the four schemes of Section~\ref{sec:four} for the
exterior Helmholtz Dirichlet problem discussed in Figures 1 and 4(c)
of~\cite{Hao13}. For convenience we repeat the details of that
problem. The curve $\gamma$ is parameterized as
\begin{equation}
r(t)=\frac{9}{20}\left(1+\frac{20}{81}\sin(5t)\right)(\cos(t),\sin(t))\,,
\quad -\pi\le t\le\pi\,.
\end{equation}
The wave number is $k=280$, which corresponds to about 48 wavelengths
across the generalized diameter of $D$. The boundary condition $g(r)$
of~(\ref{eq:HelmBV}) stem from a field excited by five point sources
with locations
\begin{equation}
r_{{\rm p}i}=a_i(\cos(b_i),\sin(b_i))\,, \quad i=1,\ldots,5\,,
\quad r_{{\rm p}i}\in D\,,
\label{eq:rpi}
\end{equation}
and with source strengths $q_i$. The values of $a_i$, $b_i$, and $q_i$
are produced by the {\sc Matlab} code (P.G. Martinsson, private
communication 2013)
\begin{verbatim}
  rand('seed',0)
  q = rand(5,1);
  a = 0.1*rand(5,1)+0.1;
  b = 2*pi*rand(5,1);
\end{verbatim}

\begin{figure}[t]
  \centering \includegraphics[height=80mm]{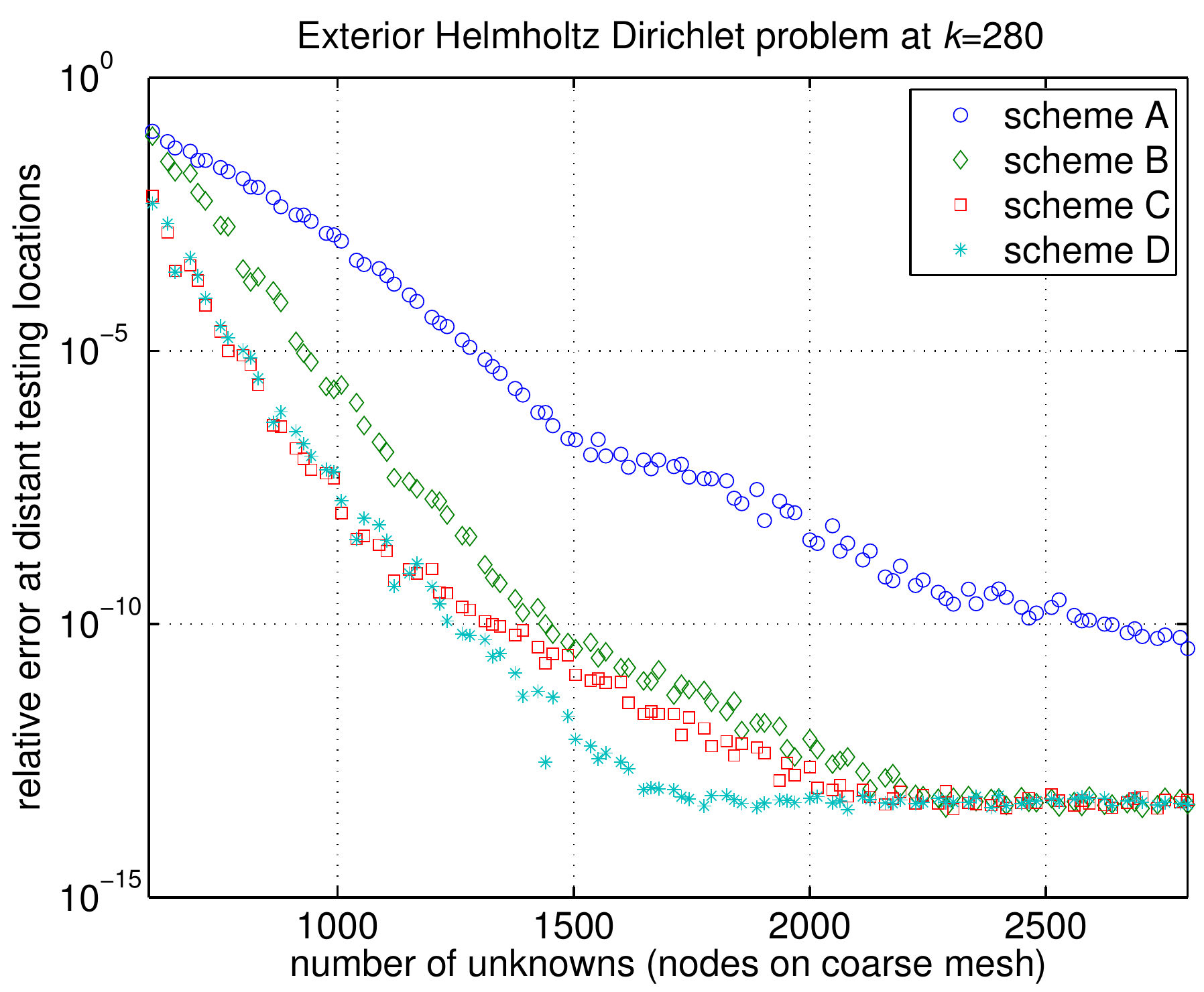}
\caption{\sf Far field tests of the schemes of Section~\ref{sec:four}.}
\label{fig:case1}
\end{figure}

We choose underlying Gauss--Legendre quadrature with $n_{\rm pt}=16$
and use the GMRES iterative solver~\cite{Saad86} for the linear
systems. The GMRES implementation involves a low-threshold stagnation
avoiding technique~\cite[Section~8]{Hels08} applicable to systems
coming from discretizations of Fredholm second kind integral
equations. The stopping criterion threshold in the (estimated)
relative residual is set to machine epsilon ($\epsilon_{\rm mach}$).

\begin{figure}[t]
\centering
\includegraphics[height=80mm]{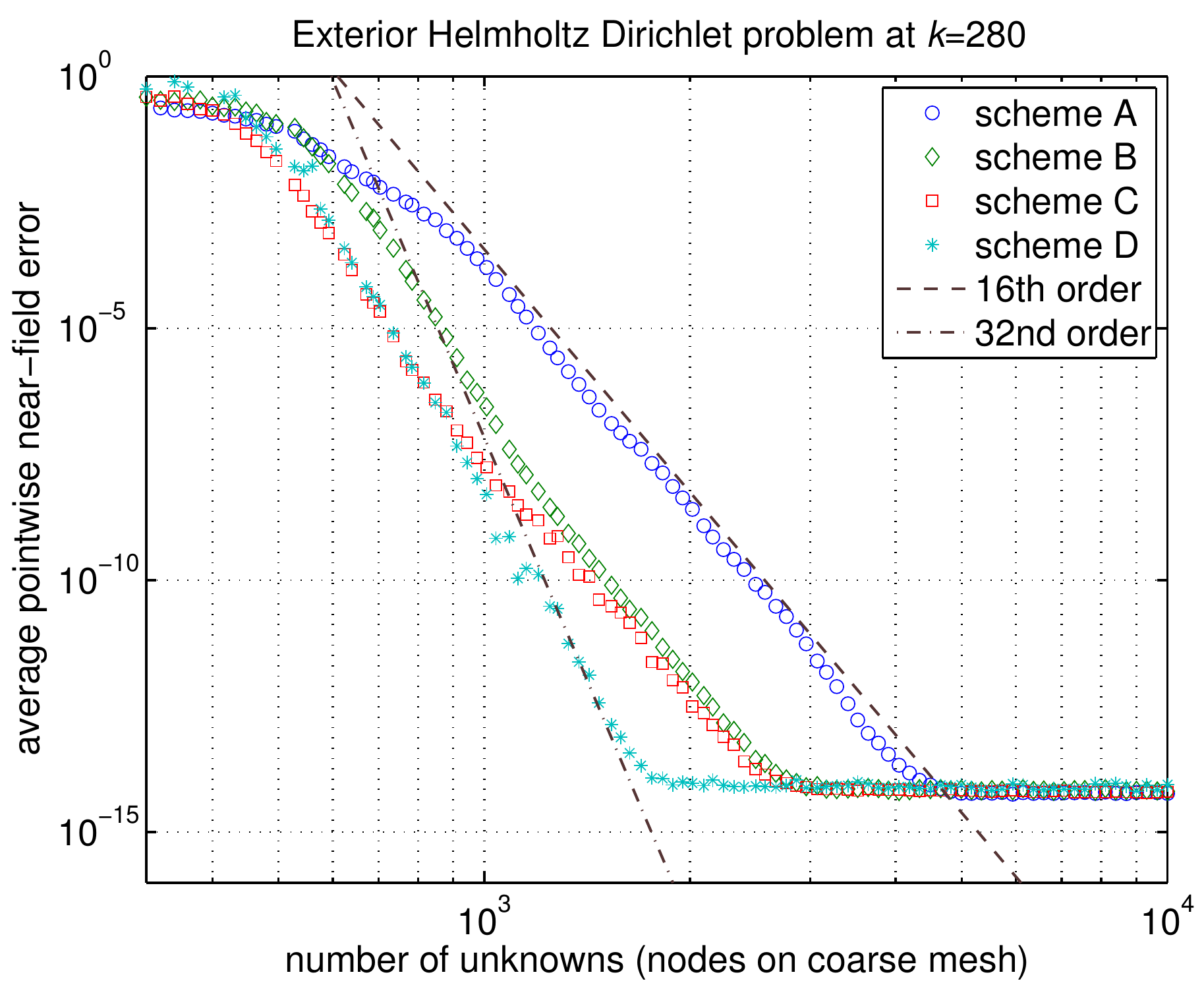}
\caption{\sf Near field tests of the schemes of Section~\ref{sec:four}.}
\label{fig:case2}
\end{figure}

\begin{figure}[t]
\centering 
\noindent\makebox[\textwidth]{
\begin{minipage}{1.1\textwidth}
\hspace{-4mm}
  \includegraphics[height=62mm]{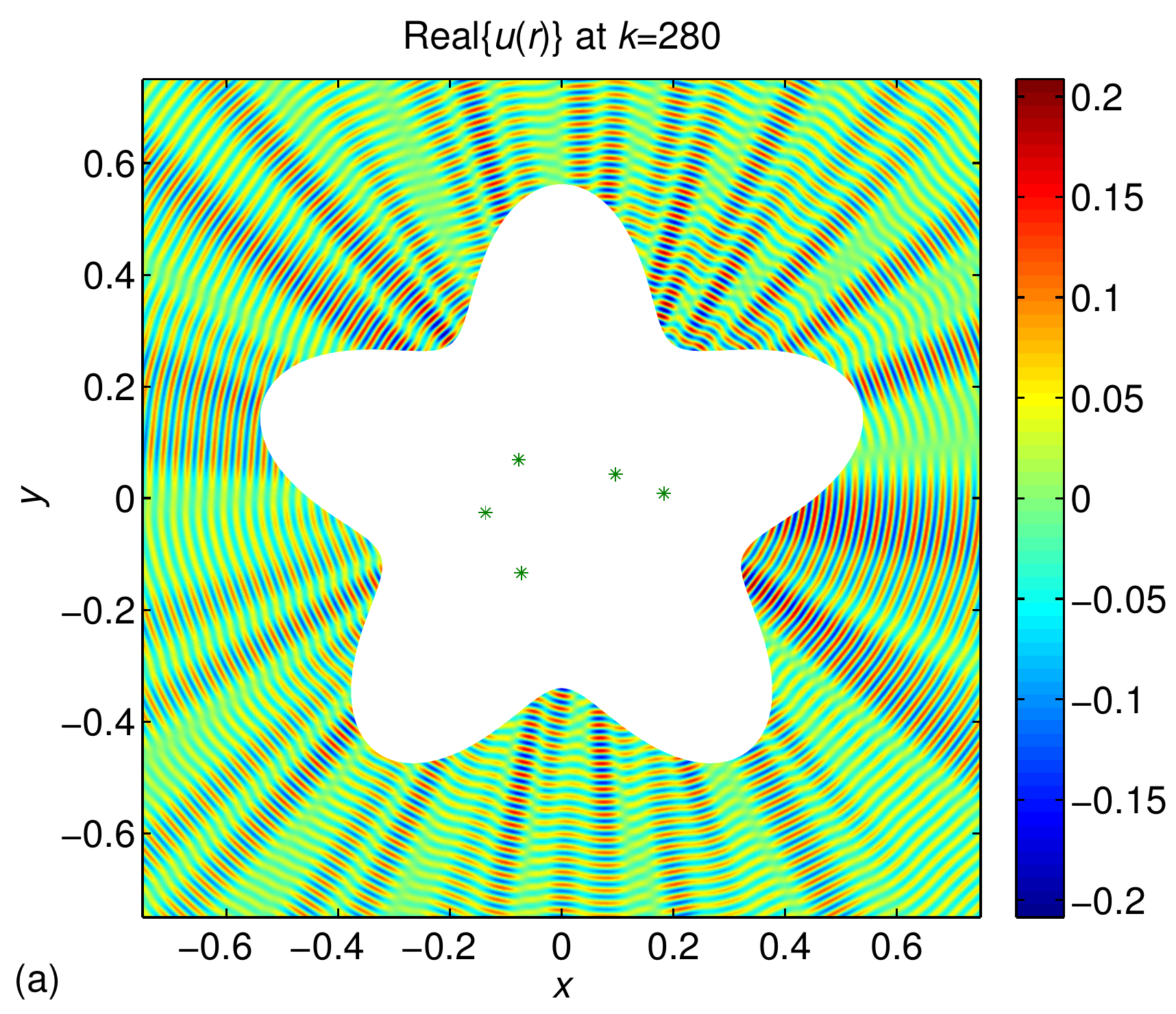}
  \includegraphics[height=62mm]{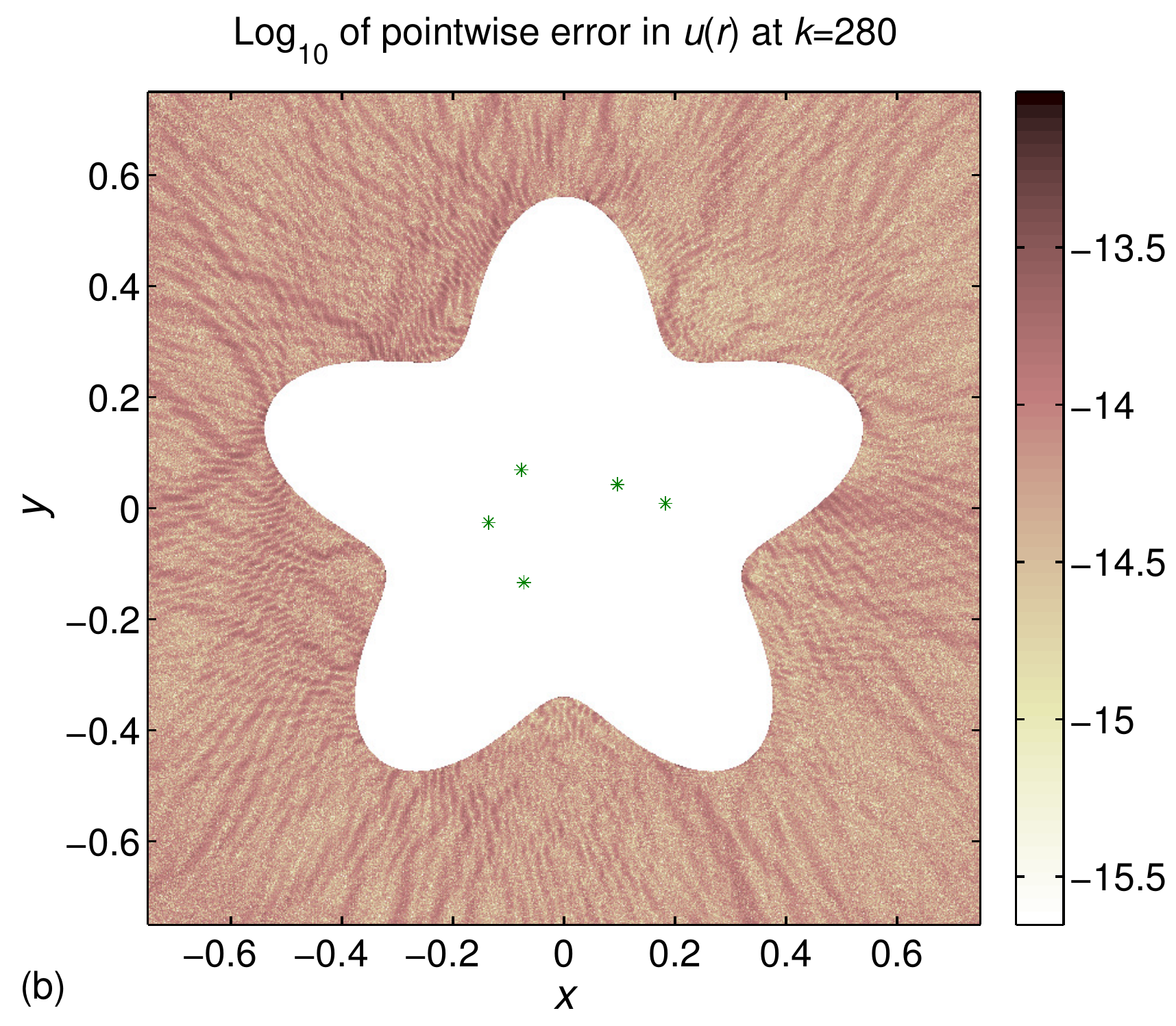}
\end{minipage}}
\caption{\sf Field and error at 347650 near-field points computed
  with scheme~C and 3904 unknowns on $\gamma$. The sources $r_{{\rm
      p}i}$ of~(\ref{eq:rpi}) that generate the boundary conditions
  are shown as green stars. (a) Real part of $u(r)$. (b) $\log_{10}$
  of pointwise error in $u(r)$ normalized with the largest value of
  $|u(r)|$ in $E$.}
\label{fig:field}
\end{figure}

Two test series are run. The first series measures the maximum
relative pointwise error at nine distant ``testing
locations''~\cite{Hao13}
\begin{equation}
r_{{\rm t}i}= 1.25(\cos(2\pi(i-1)/9),\sin(2\pi(i-1)/9))\,,
\quad i=1,\ldots,9\,.
\end{equation}
The second series measures the average pointwise error, normalized
with the largest value of $|u(r)|$, in a near-field zone. This zone is
taken as the intersection of $E$ and the square $x,y\in[-0.75,0.75]$,
where a Cartesian grid of $200\times 200$ field points $r$ is created.
Points in $D$ are then excluded, leaving a number of 28,460 points
where $u(r)$ is evaluated. Test results for the far field are
presented in Figure~\ref{fig:case1} and for the near field in
Figure~\ref{fig:case2}.

Figure~\ref{fig:case1}, where the $x$-axis has linear scaling as to
facilitate comparison with Figure~4(c) of~\cite{Hao13}, shows that
fine grid resolution of $M(r,r')$ for $r'$ close to $r$ (scheme~B)
substantially improves on pure coarse grid resolution (scheme~A). The
number of unknowns needed to resolve $\rho(r)$, for a given accuracy
in the far field $u(r)$, is roughly cut in half. The benefit of using
unit speed parameterization (scheme~C) is clearly visible too.
Extended interpolation (scheme~D) is only worthwhile when the highest
achievable accuracy is of interest. A comparison between our scheme B
and the 10th order accurate ``Modified Gaussian'' scheme
of~\cite{Hao13} reveals that scheme~B converges faster, as expected.
For example, with 1600 discretization points on $\gamma$ the gain in
accuracy is around three digits. The achievable accuracy of those of
our schemes that have saturated is on par with that of Kress' global
explicit-split scheme~\cite{Kress91}, which is the most accurate of
the schemes tested in~\cite{Hao13} and which also exhibits a more
rapid convergence than any of our schemes in this example.

The results of the near field tests in Figure~\ref{fig:case2} are
similar to those of the far field tests. The 16th order convergence of
scheme A is apparent thanks to the logarithmic scaling of the
$x$-axis. Schemes B and C also exhibit asymptotic 16th order
convergence while scheme~D converges more rapidly and with no evident
asymptotics visible. Further experiments (not shown) indicate that the
choice $n_{\rm s}=4$ in scheme~D is optimal in the sense that larger
$n_{\rm s}$ do not improve the convergence.

The actual field $u(r)$ and an example of a near-field error plot
produced by scheme~C are shown in Figure~\ref{fig:field}(a)
and~\ref{fig:field}(b). A Cartesian grid of $700\times 700$ field
points is used to produce these images. One can see, in
Figure~\ref{fig:field}(b), that the field error is uniformly small all
the way up to the boundary $\gamma$ although the convergence close to
concave parts of $\gamma$ is, in fact, somewhat delayed compared to
that in the remainder of $E$. The field point that happens to lie
closest to $\gamma$ in this example is only separated from $\gamma$ by
a distance of $3\cdot 10^{-6}$. Further tests (not shown) indicate
that our product integration scheme can evaluate $u(r)$ at $r\in E$
arbitrarily close to $\gamma$ without the error deteriorating except
for when $r$ lies extremely close to an endpoint of a panel. Then
cancellation occurs in the quantities {\tt p(1)} and {\tt p1} of
Appendix~B. In this case a procedure involving temporary panels
mergers and splits can restore the
accuracy~\cite[Section~5.5]{Hels08}.

The achievable pointwise precision for $u(r)$ ranges between 13 and 15
digits in all our resolved examples. This is comparable to the
precision obtained by Barnett in a similar example, exterior to a
domain 12 wavelengths in diameter, using Kress' global explicit-split
scheme with a QBX post-processor~\cite[Section 4.1]{Barn14}. Given
that the condition numbers of the main system matrices ${\bf I}+{\bf
  M}_{\gamma}$ in our schemes are very low, less than eight for
$k=280$, one could speculate that it is possible to construct even
more accurate schemes.

We end with some comments on the convergence of the GMRES iterative
solver and on the influence of the coupling parameter $\eta$,
discussed in Remark~\ref{rmk:eta}. With $\eta=k/2$ (our preferred
choice) and the stopping criterion threshold set to $\epsilon_{\rm
  mach}$, GMRES converges in 51 iterations for all schemes and most
resolutions in Figures~\ref{fig:case1} and~\ref{fig:case2}. Severely
underresolved systems require more iterations, though. With the choice
$\eta=k$, made in~\cite[Equation~(7.7)]{Hao13}, the typical number of
iterations required for full convergence rises to 60. The choice
$\eta=-k$, made in the work on direct solvers for 3D scattering
problems~\cite[Equation~(2.3)]{Brem14}, gives a particular slow
convergence in GMRES -- as also noted in~\cite[Remark~6.1]{Brem14}. A
number of 373 iterations is needed in our experiments. In order to
compare our GMRES convergence results to those
of~\cite[Table~2]{Hao13}, we increase the stopping criterion threshold
to $10^{-12}$ and lower the wavenumber to $k=2.8$. Then 13 iterations
are needed. This is marginally better than the 14 iterations reported
for all competitive schemes in~\cite[Table~2]{Hao13}.

\section{Discussion}

Initially, the thought of implementing an explicit kernel-split
panel-based Nyström discretization scheme for planar Helmholtz
boundary value problems may seem daunting because of the hefty series
representations of Hankel functions that one can encounter when
searching the special functions literature for information about the
kernel singularities of the single- and double layer operators $S_k$
and $K_k$. But, when the dust has settled, the kernel splits can be
written out in amazingly simple forms which, essentially, only involve
functions that need to be evaluated also in split-free schemes. With
access to efficient product integration techniques for logarithmic-
and Cauchy-type singular kernels, the implementation is very
straightforward.

The present paper shows that, for planar problems, the explicit
kernel-split philosophy with quadratures computed on the fly offers
rapid and stable convergence for linear systems as well as for far
field and near field solutions. It, further, allows for the
construction of schemes where mixes of discretizations on coarse and
fine meshes enhance the performance. 
 
In three dimensions the situation is more involved. For one thing,
product integration for singular kernels is problematic. While it is
possible to carry the techniques of the present work over to
axisymmetric Helmholtz problems (a modification of scheme~B is used
in~\cite{Hels14}), it is quite likely that schemes relying on adaptive
and precomputed quadratures are better suited to solving fully three
dimensional problems~\cite{Brem14}.

\bigskip\noindent 
{\bf Acknowledgements.} Johan Helsing wishes to thank the Banff
International Research Station, Alberta, Canada (BIRS) and its staff
for providing a creative atmosphere at the workshop ``Integral
Equations Methods: Fast Algorithms and Applications'' in December 2013
where parts of this work were finalized. The work was supported by the
Swedish Research Council under grant 621-2011-5516.

\newpage
\section*{Appendix A: code for $\mathfrak{W}_{\rm L}$}

The following {\sc Matlab} function returns the $n_{\rm pt}\times
n_{\rm pt}$ matrix $\mathfrak{W}_{\rm L}$, needed for the construction
of the product integration weight corrections $w_{{\rm L}j}^{\rm
  corr}(r_i)$ of~(\ref{eq:Acorr}):
\begin{verbatim}
  function WfrakL=WfrakLinit(trans,scale,tfrak,npt)
  A=fliplr(vander(tfrak));
  tt=trans+scale*tfrak;
  Q=zeros(npt);
  p=zeros(1,npt+1);
  c=(1-(-1).^(1:npt))./(1:npt);
  for m=1:npt
    p(1)=log(abs((1-tt(m))/(1+tt(m))));
    p1=log(abs(1-tt(m)^2));
    for k=1:npt
      p(k+1)=tt(m)*p(k)+c(k);
    end
    Q(m,1:2:npt-1)=p1-p(2:2:npt);
    Q(m,2:2:npt)=p(1)-p(3:2:npt+1);
    Q(m,:)=Q(m,:)./(1:npt);  
  end
  WfrakL=Q/A;
\end{verbatim}
The choice of input arguments {\tt trans=0} and {\tt scale=1} gives
$\mathfrak{W}_{\rm L}$ for $r_i$ and $r_j$ on the same quadrature
panel $\gamma_p$. The choice {\tt trans=$\pm$2} and {\tt scale=1}
gives $\mathfrak{W}_{\rm L}$ for $r_i$ on a neighboring panel
$\gamma_{p\pm 1}$, assuming it is equal in parameter length. The input
argument {\tt tfrak} is a column vector whose entries contain the
canonical nodes $\mathfrak{t}_i$, $i=1,\ldots,n_{\rm pt}$, and {\tt
  npt} corresponds to $n_{\rm pt}$.

Note that special-purpose quadrature only must be activated when $r_i$
and $r_j$ are close to each other, see Section~\ref{sec:active}, and
that no closeness check is included in ${\tt WfrakLinit}$.
Furthermore, a marginal improvement in accuracy can result from
running the recursion for {\tt p} backwards in certain situations,
see~\cite[Section 6]{Hels08}.

\newpage
\section*{Appendix B: code for $w_{{\rm L}j}^{\rm corr}(r)$
  and $w_{{\rm C}j}^{\rm cmp}(r)$}

The following {\sc Matlab} function returns the weight corrections
$w_{{\rm L}j}^{\rm corr}(r)$ and the compensation weights $w_{{\rm
    C}j}^{\rm cmp}(r)$, $j=1,\ldots,n_{\rm pt}$, needed
in~(\ref{eq:fix2}) and~(\ref{eq:fix3}) when $r\in E$:
\begin{verbatim}
  function [wcorrL,wcmpC]=wLCinit(ra,rb,r,rj,nuj,rpwj,npt)
  dr=(rb-ra)/2;
  rtr=(r-(rb+ra)/2)/dr;
  rjtr=(rj-(rb+ra)/2)/dr;
  A=fliplr(vander(rjtr)).';
  p=zeros(npt+1,1);
  q=zeros(npt,1);
  c=(1-(-1).^(1:npt))./(1:npt);
  p(1)=log(1-rtr)-log(-1-rtr);
  p1=log(1-rtr)+log(-1-rtr);    
  if imag(rtr)>0 && abs(real(rtr))<1
    p(1)=p(1)-2i*pi;
    p1=p1+2i*pi;        
  end
  for k=1:npt
    p(k+1)=rtr*p(k)+c(k);
  end
  q(1:2:npt-1)=p1-p(2:2:npt);
  q(2:2:npt)=p(1)-p(3:2:npt+1);
  q=q./(1:npt)';
  wcorrL=imag(A\q*dr.*conj(nuj))./abs(rpwj)-log(abs((rj-r)/dr));
  wcmpC=imag(A\p(1:npt)-rpwj./(rj-r));
\end{verbatim}
This function relies on complex arithmetic and takes, as input, points
and vectors in $\mathbb{R}^2$ represented as points in $\mathbb{C}$.
Otherwise the notation follows Section~\ref{sec:prod}: input
parameters {\tt ra} and {\tt rb} correspond to $r(t_a)$ and $r(t_b)$;
{\tt r} is the target point $r\in E$; and {\tt rj}, {\tt nuj}, and
{\tt rpwj} are column vector whose entries contain the points $r_j$,
the exterior unit normals $\nu$ at $r_j$, and the weighted velocity
function $\dot{r}_jw_j$, $j=1,\ldots,n_{\rm pt}$.

\end{document}